\documentclass[UTF8]{article}
\usepackage{graphicx} 
\usepackage{amsmath} 
\usepackage{amssymb} 
\usepackage{setspace}

\title{Existence, Nonexistence, and Symmetry of Positive Solutions for Fractional Laplacian Problems}
\author{Lu Haipeng, Yu Mei}

\begin{document}

      \maketitle

\begin{abstract}
    This paper studies the properties of solutions to a class of elliptic and parabolic problems involving the fractional Laplacian. By applying the mountain pass theorem, we prove the existence of bounded classical positive solutions in the subcritical regime. Moreover, using the method of moving planes, we establish that these solutions are symmetric or monotone in the first variable. In contrast, we show that no such solutions exist in the supercritical or negative exponent cases. An analysis of the asymptotic behavior of solutions at infinity provides further insight into their profiles, which supports applications to real-world problems. The approaches developed in this work can also be extended to a wider range of nonlocal elliptic and parabolic equations, including those with more general operators and nonlinearities.

\end{abstract}

\section{Introduction}\label{sec1}

\quad \quad We here analyze the properties of solutions to the fractional elliptic equation
\[(-\Delta)^{s}u(x)=x_{1}u^{p}(x)+f(u(x)),\quad x\in\Omega, \tag{1.1}\]
and the fractional parabolic equation
\[\frac{\partial u}{\partial t}(x,t)+(-\Delta)^{s}u(x,t)=x_{1}u^{p}(x,t)+f(u(x,t) ),\quad (x,t)\in\Omega\times\mathbb{R}, \tag{1.2}\]
where \(\Omega\) is either a bounded domain or the half-space \(\mathbb{R}^{n}_{1,+}:=\{x\in\mathbb{R}^{n} \mid x_{1}>0\}.\) Here the fractional Laplacian operator $(-\Delta)^s$ is nonlocal and defined via the singular integral
\[(-\Delta)^{s}u(x):=C_{n,s}\,\mathrm{P.V.}\int_{\mathbb{R}^{n}}\frac{u(x)-u(y)}{|x -y|^{n+2s}}\mathrm{d}y, \tag{1.3}\]
where \( s \in (0,1)\) and \(\mathrm{P.V.}\) denotes the Cauchy principal value. The integral in (1.3) is well-defined for functions \(u \in L_{2s} \cap C^{1,1}_{\mathrm{loc}},\) with the space \(L_{2s}\) given by
\[L_{2s}:=\left\{u\in L^{1}_{\mathrm{loc}}(\mathbb{R}^{n})\;\middle|\;\int_{\mathbb{R }^{n}}\frac{|u(x)|}{1+|x|^{n+2s}}\mathrm{d}x<+\infty\right\}.\]

Due to their capacity to capture long-range interactions, the nonlocal operators, especially fractional Laplacian, are often more suitable than the standard Laplacian or variable-coefficient elliptic operators for modeling various phenomena. The key advantage of the fractional Laplacian lies in its ability to accurately describe complex systems that exhibit long-range interactions, nonlocal correlations, and fractal structures. As a result, such equations are increasingly applied in diverse areas such as quasi-geostrophic flows \cite{1}, activator-inhibitor systems \cite{2}, nonlocal diffusion processes \cite{3,4}, and fractional porous media flows \cite{5}. It is crucial for characterizing their asymptotic behavior as spatial variables approach infinity to understand the monotonicity and symmetry of solutions, thereby enhancing the applicability of the models to real-world problems. Moreover, proving the existence of solutions is fundamental for validating the underlying mathematical frameworks.

Research on the properties of solutions to fractional equations has developed into a relatively comprehensive theoretical framework. In the context of elliptic equations, \cite{6} use the method of moving planes to establish radial symmetry criteria for \((-\Delta)^{\alpha/2}u=u^{p},\) highlighting the critical role of the exponent \(p=\frac{n+\alpha}{n-\alpha}.\) They also show that solutions to the Dirichlet problem in a half-space must vanish in the subcritical case. Furthermore, \cite{7} prove that solutions within a spherical domain exhibit radial symmetry about the origin (under the assumption that \(f\) is Lipschitz continuous), while traveling wave solutions in full space are strictly monotonic along the \(x_1\)-direction (requiring \(f\) to be piecewise non-increasing). \cite{8} study the concave-convex problem involving the fractional Laplacian
\[\begin{cases}
(-\Delta)^{\frac{\alpha}{2}}u=\lambda u^{q}+u^{p}, & u>0 \quad \text{in } \Omega,\\ 
u=0 & \text{on } \partial\Omega,
\end{cases}\tag{1.4}\]
with \(0<\alpha<2, 0<q<1<p<\frac{N+\alpha}{N-\alpha}, N>\alpha, \lambda>0\), and \(\Omega\subset \mathbb{R}^{N}\) is a smooth bounded domain. By employing methods such as local extension, variational frameworks, and the sub- and super-solution technique, they conclude that the number of solutions to (1.4) depends on \(\lambda\), and that solutions are bounded when \(\alpha>1.\) \cite{9} investigate the Dirichlet problem for a semilinear elliptic system involving the fractional Laplacian
\[\begin{cases}
(-\Delta)^{\frac{\alpha}{2}}u(x)=v^{q}(x), & x\in\mathbb{R}_{+}^{n},\\ 
(-\Delta)^{\frac{\alpha}{2}}v(x)=u^{p}(x), & x\in\mathbb{R}_{+}^{n},\\ 
u(x)=v(x)=0, & x\notin\mathbb{R}_{+}^{n},
\end{cases}\tag{1.5}\]
for \(\alpha \in (0,2)\) and \(\frac{n}{n-\alpha}<p,q\leq\frac{n+\alpha}{n-\alpha},\) with \(u,v\in L_{\rm loc}^{\frac{n(p-1)}{\alpha}}(\mathbb{R}_{+}^{n})\times L_{\rm loc}^ {\frac{n(q-1)}{\alpha}}(\mathbb{R}_{+}^{n}).\) By converting the system into an integral form and applying the method of moving planes in integral form among other techniques, they show that any pair of nonnegative solutions \((u,v)\) of (1.5) must be trivial.

In the parabolic setting, \cite{10} analyze the \(\omega\)-limit set and binary reduction properties of solutions in the unit ball (either trivial or radially monotone about the origin) and demonstrate that solutions in full space exhibit spatially localized symmetry when \(f_{u}(t,0)<-\sigma.\) \cite{11} make a significant breakthrough by discovering translation symmetry phenomena (i.e., the existence of a direction such that the solution is symmetric with respect to it). In the study of new operators, \cite{12} first confirm that solutions to the full fractional heat operator equation \((\partial_{t}-\Delta)^{s}u=x_{1}u^{p}\) are strongly constrained along the x-direction (leading to non-existence when \(p>1\)). Their analysis of mixed operators \cite{13} reveal that divergence in the spatial coefficient \(a(x)\) can preclude the existence of solutions, while directional monotonicity is governed by the coupling condition \(a(x)f(u).\) These results establish a universal connection between the properties of nonlocal operators and the behavior of solutions.

The aforementioned contributions are vital for advancing the theory of fractional equations and enhancing their applicability to real-world problems. However, the nonlinear terms in these studies are often restrictive, limiting their ability to model more complex scenarios. To address this gap, this paper investigates the properties of solutions to fractional Laplacian elliptic and parabolic problems with a specific nonlinearity of the form \(x_{1}u^{p}+f(u).\)

The rest of this paper is structured as follows. Section 2 is devoted to the fractional Laplacian elliptic problems. We first prove the existence of solutions in the subcritical case via the mountain pass theorem and energy functionals. Then we apply the method of moving planes to establish monotonicity and symmetry (Theorems 2.1.1 and 2.1.2). Finally, we prove the nonexistence of nontrivial classical positive solutions in the supercritical and negative exponent cases using Pohozhaev-type identities and boundary analysis (Theorems 2.2.1 and 2.3.1). Section 3 addresses the fractional parabolic problems. We employ tools such as Galerkin approximations, energy estimates, and boundary growth estimates to establish analogous existence, monotonicity, symmetry, and nonexistence results, corresponding to Theorems 3.1.1, 3.1.2, 3.2.1, and 3.3.1. Finally, in Section 4, we provide a concluding summary of the main findings, discuss the implications of our work, and suggest potential directions for future research.

\section{Investigation on Properties of Solutions to Fractional Laplacian Elliptic Problems}\label{sec2}

\quad \quad This section mainly presents the existence, monotonicity, symmetry, and other related properties of solutions to the fractional Laplacian elliptic problem
\[\begin{cases}
(-\Delta)^s u(x) = x_1 u^p(x) + f(u(x)), & u \in H^s(\Omega), u > 0, \quad x \in \Omega \subset \mathbb{R}^n \\
u(x) \equiv 0, & x \in \Omega^c
\end{cases} \tag{2.1}\]
where \(\Omega = B_1(0)\) or \(\mathbb{R}_{1,+}^n\). Here, \(f: H^s \to \mathbb{R}\) is a continuous, integrable functional, the fractional order  \(s \in (0,1)\), and the space dimension \(n \geq 2\).

Regarding the exponent \(p\), we will discuss the following cases:
\begin{itemize}
    \item Case I (Subcritical case): \(0 < p < \frac{n+2s}{n-2s}\);
    \item Case II (Supercritical case): \(p \geq \frac{n+2s}{n-2s}\);
    \item Case III(Negative exponent case): \(p < 0\).
\end{itemize}

\subsection{Properties of Solutions in the Subcritical Case}

\quad \quad This subsection shows the properties of solutions to problem (2.1) when \(0 < p < \frac{n+2s}{n-2s}\) and \(\Omega = B_1(0)\) or \(\mathbb{R}_1^n\). First, we lead in the following lemmas, which are largely related with our conclusions:

\textbf{Lemma 2.1.1.} [Narrow Region Principle \cite{7}]  Suppose there exists a point \( x^0 \in \Sigma_{\lambda},\)  such that \( w_\lambda(x^0) = \min_{\Sigma_{\lambda}} w_\lambda(x),\)  then 
\[(-\Delta)^s w_\lambda(x^0) \leq \frac{C_1}{[d(x^0,T_\lambda)]^{2s}} w_\lambda(x^0). \tag{2.2}\]
where \(d(x^0,T_\lambda)\) is the distance from point \(x^0\) to the plane \(T_\lambda\).

\textbf{Lemma 2.1.2.} [Maximum Principle for Unbounded Domains \cite{14,15}] Let \( \Omega \subset \mathbb{R}^n \) be an open set that may be unbounded or disconnected, satisfying for any \(x \in \mathbb{R}^n,\)
\[\lim_{R \to +\infty} \frac{|B_R(x) \cap \Omega^c|}{|B_R(x)|} \geq c_0 > 0.\]
Assume that a function \( u \in C^{1,1}_{loc}(\Omega) \cap L^{2s}(\Omega) \) upper semicontinuous on \(\Omega\), bounded from above, and satisfying
\[\begin{cases}
(-\Delta)^s u(x) \leq 0, & x \in \Omega \text{ and } u(x) > 0, \\
u(x) \leq 0, & x \in \Omega^c,
\end{cases}\]
then \( u(x) \leq 0 \) on \(\Omega\).

Based on the above lemmas, we have the following results:

\textbf{Theorem 2.1.1.} 

(i) Suppose \( 0 < p < \frac{n+2s}{n-2s} \), then the problem (2.1) admits a bounded classical positive solution \( u \in  C^{2s} \left(B_1(0)\right) \).

(ii) Let \( u \in C^{2s} \left(B_1(0) \right) \) be a solution to problem (2.1). Assume \( f \) is Lipschitz continuous, and for some \( q > p \) and a constant \( C > 0 \) depending only on \( n, s, q \), we have
\[|f(u)| \leq C|u|^q.\]
Then \( u \) is symmetric with respect to the hyperplane \( \{x \in \mathbb{R}^n | x_1 = 0\} \), and monotone decreasing in \( x_1 \) when \( x_1 > 0 \).

\textbf{Proof.} (i) First, we construct an energy functional \( J(u) \) of (2.1):
\[J(u) = \int_{B_1(0)} \left( \frac{1}{2} |(-\Delta)^{\frac{s}{2}} u|^2 - F(u) - \frac{x_1}{p+1} u^{p+1} \right) dx, u \in H^s(B_1(0)),\tag{2.3}\]
where \( F(u) = \int_0^u f(v) dv \). Next, verify that \( J(u) \) is bounded below in \( u \in H^s(B_1(0)) \). Owing to the fractional Sobolev inequality \cite{16} and the Gagliardo-Nirenberg interpolation inequality \cite{17}, we obtain
\[\|u\|_{L^{p+1}(B_1(0))} \leq C \|u\|_{H^s(B_1(0))}. \tag{2.4}\]
From the growth condition of \( f \), we have:
\[|F(u)| \leq C(1 + |u|^{q+1}). \tag{2.5}\]
The equations (2.3)-(2.5) yield that
\[J(u) \geq \frac{1}{2} \|(-\Delta)^{\frac{s}{2}} u\|^2_{L^2(B_1(0))} - C \left( 1 + \|u\|^{q+1}_{L^{p+1}(B_1(0))} \right) - \|u\|^{p+1}_{L^{p+1}(B_1(0))}. \]
Using appropriate Young's inequality \cite{18} and the Sobolev embedding theorem \cite{19}, we can prove that \(J(u)\) is bounded below in \(H^s(B_1(0))\).

Then we need to prove that \(J(u)\) is lower semicontinuous under weak topology. Suppose \(u_k \rightharpoonup u\) weakly in \(H^s(B_1(0))\), then
\[\liminf_{k \to \infty} \|(-\Delta)^{\frac{s}{2}} u_k\|_{L^2(B_1(0))}^2 \geq \|(-\Delta)^{\frac{s}{2}} u\|_{L^2(B_1(0))}^2. \tag{2.6}\]
For the nonlinear terms \(F(u)\) and \(u^{p+1}\), we can prove that
\[\liminf_{k \to \infty} \int_{B_1(0)} F(u_k) dx \geq \int_{B_1(0)} F(u) dx, \tag{2.7}\]
\[\liminf_{k \to \infty} \int_{B_1(0)} u_k^{p+1} dx \geq \int_{B_1(0)} u^{p+1} dx. \tag{2.8}\]
by Fatou's lemma \cite{18}.

From (2.6)-(2.8), it follows that the energy functional \(J(u)\) is lower semicontinuous under weak topology. Therefore, we can construct a minimizing sequence \(u_k \subset H^s(B_1(0))\) such that \(J(u_k) \to \inf\limits_{u \in H^s(B_1(0))} J(u)\). By Banach-Alaoglu theorem, we can extract a weakly convergent subsequence (we still denote as \(u_k\)) such that \(u_k \rightharpoonup u\) in \(H^s(B_1(0))\). Since \(J(u)\) is lower semicontinuous, we have:
\[J(u) \leq \liminf_{k \to \infty} J(u_k) = \inf_{u \in H^s(B_1(0))} J(u).\]
Therefore, \(u\) is a critical point of \(J(u)\).

Finally, we shall prove that \(u\) is a classical solution in \( B_1(0)\). By the fractional Sobolev embedding theorem \cite{19} and \(s \in (0,1)\), we have \(s < \frac{n}{2}\), so the embedding yields \(u \in L^q(B_1(0))\), where \(q = \frac{2n}{n-2s}\). Through De Giorgi-Nash-Moser estimates \cite{20}, we can improve the regularity of the solution to \(u \in C^{2s}(B_1(0))\). Thus \(u \in C^{2s}(B_1(0))\) is a classical solution.

(ii) The proof follows a similar line of reasoning as presented by \cite{7} and is included here for completeness. Denote the following notations:
\begin{itemize}
    \item \(T_\lambda := \{x \in \mathbb{R}^n | x_1 = \lambda\}\) is the moving plane;
    \item \(\Sigma_\lambda := \{x \in \mathbb{R}^n | x_1 < \lambda\}\) is the region to the left of the moving plane;
    \item \(x^\lambda := (2\lambda - x_1, x_2, \cdots, x_n)\) is the reflection point with respect to the hyperplane \(T_\lambda\);
    \item \( u_{\lambda}(x) := u(x^{\lambda}) \);
    \item \( w_{\lambda}(x) := u_{\lambda}(x) - u(x) \);
    \item \(\Omega_{\lambda} := \Sigma_{\lambda} \cap B_{1}(0)\).
\end{itemize}

We first verify that for \(\lambda\) sufficiently close to \(-1\), 
\[w_{\lambda}(x) \geq 0, \forall x \in \Omega_{\lambda}. \tag{2.9}\]

Suppose (2.9) does not hold, then there exists a point \(x^{0} \in \Omega_{\lambda}\) such that:
\[w_{\lambda}(x^{0}) = \min_{x \in \Omega_{\lambda}} w_{\lambda}(x) < 0.\]
From (2.1) we have
\[(-\Delta)^{s}u(x) = x_{1}u^{p}(x) + f(u(x)), \tag{2.10}\]
\[(-\Delta)^{s}u_{\lambda}(x) = (x^{\lambda})_{1}u_{\lambda}^{p}(x) + f(u_{\lambda}(x)). \tag{2.11}\]
Subtracting (2.10) from (2.11) yields
\begin{align}
(-\Delta)^{s}w_{\lambda}(x) &= (-\Delta)^{s}u_{\lambda}(x) - (-\Delta)^{s}u(x) \nonumber \\
&= [(x^{\lambda})_{1}u_{\lambda}^{p}(x) + f(u_{\lambda}(x))] - [x_{1}u^{p}(x) + f(u(x))] \nonumber \\
&= [(x^{\lambda})_{1}u_{\lambda}^{p}(x) - x_{1}u^{p}(x)] + [f(u_{\lambda}(x)) - f(u(x))]. \tag{2.12}
\end{align}
Since
\begin{align}
(x^{\lambda})_1u_{\lambda}^{p}(x) - x_1u^{p}(x) &= (x^{\lambda})_1u_{\lambda}^{p}(x) - x_1u_{\lambda}^{p}(x) + x_1u_{\lambda}^{p}(x) - x_1u^{p}(x) \nonumber \\
&= [(x^{\lambda})_1 - x_1]u_{\lambda}^{p}(x) + x_1[u_{\lambda}^{p}(x) - u^{p}(x)] \nonumber \\
&\geq x_1[u_{\lambda}^{p}(x) - u^{p}(x)] \nonumber \\
&= x_1 \cdot p \xi_\lambda^{p-1}(x)[u_{\lambda}(x) - u(x)] \nonumber \\
&= x_1 \cdot p \xi_\lambda^{p-1}(x) w_{\lambda}(x), \tag{2.13}
\end{align}
where \(\xi_{\lambda}(x)\) lies between \(u_{\lambda}(x)\) and \(u(x)\) for each \(x\). Similarly, 

\begin{align}
    f(u_{\lambda}(x)) - f(u(x)) &= \frac{f(u_{\lambda}(x)) - f(u(x))}{u_{\lambda}(x) - u(x)} \cdot [u_{\lambda}(x) - u(x)] \nonumber \\ &= \frac{f(u_{\lambda}(x)) - f(u(x))}{u_{\lambda}(x) - u(x)} \cdot w_{\lambda}(x). \tag{2.14}
\end{align}

Substituting (2.13) and (2.14) into (2.12), we get:
\begin{align}
(-\Delta)^s w_\lambda(x) &\geq \left[ x_1 \cdot p_{\xi_{\lambda}}^{p-1}(x) + \frac{f(u_\lambda(x)) - f(u(x))}{u_\lambda(x) - u(x)} \right] w_\lambda(x) \nonumber \\
&:= c_\lambda(x) w_\lambda(x). \nonumber
\end{align}
By the boundedness of \(\Omega_\lambda\) and the Lipschitz continuity of \(f\), we have \(c_\lambda(x)\) is a bounded function. Then \(w_\lambda(x)\) satisfies
\[\begin{cases}
(-\Delta)^s w_\lambda(x) \geq c_\lambda(x)w_\lambda(x), & x \in \Omega_\lambda, \\
w_\lambda(x) = -w_\lambda(x^\lambda), & x \in \Sigma_\lambda, \\
w_\lambda(x) \geq 0, & x \in \Sigma_\lambda \setminus \Omega_\lambda.
\end{cases}\]
By (2.2), we get
\[\frac{C_1}{[d(x^0,T_\lambda)]^{2s}} w_\lambda(x^0) \geq c_\lambda(x^0) w_\lambda(x^0).\]
According to the assumption \(w_\lambda(x^0) < 0\), we have
\[\frac{C_1}{[d(x^0,T_\lambda)]^{2s}} \leq c_\lambda(x^0).\]
If \(\lambda\) is sufficiently close to \(-1\), then \(\Omega_\lambda\) is a narrow region, so the distance \(d(x^0,T_\lambda)\) can be made sufficiently small, contradicting the boundedness of \(c_\lambda(x)\). Thus (2.9) holds.

Finally, we move the plane \(T_\lambda\) to the right to its limiting position while maintaining \(w_\lambda \geq 0\). Denote
\[\lambda_0 = \sup \left\{ \lambda < 0 \mid w_\mu(x) \geq 0, x \in \Omega_{\mu}, \mu \leq \lambda \right\}.\]
Clearly \(\lambda_0 \leq 0\), we shall prove \(\lambda_0 = 0\) by contradiction.

Suppose \(\lambda_0 < 0\). By the strong maximum principle, we have
\[w_{\lambda_0}(x) > 0, \forall x \in \Omega_{\lambda_0}.\]
Move the plane back by an arbitrarily small distance \(\delta > 0\), then:
\[w_{\lambda_0} \geq c_0 > 0, \forall x \in \Omega_{\lambda_0-\delta}.\]
By the continuity of \(w_\lambda\) with respect to \(\lambda\), there exists \(\epsilon > 0\) such that for any \(\lambda \in [\lambda_0, \lambda_0 + \epsilon)\), we have
\[w_\lambda(x) \geq 0, x \in \Omega_{\lambda_0-\delta} \tag{2.15}\]
holds. Since \(\Omega_{\lambda_0-\delta}\) is a narrow region, by the narrow region principle, it holds that
\[w_\lambda(x) \geq 0, x \in \Omega_\lambda \setminus \Omega_{\lambda_0-\delta}. \tag{2.16}\]

Combining (2.15) and (2.16) shows that for any \(\lambda \in [\lambda_0, \lambda_0 + \epsilon)\), we have:
\[w_\lambda(x) \geq 0, x \in \Omega_\lambda,\]
which contradicts the definition of \(\lambda_0\), therefore \(\lambda_0 = 0\). Consequently, \(u(x)\) is symmetric with respect to the hyperplane \(\{x \in \mathbb{R}^n | x_1 = 0\}\), and monotone decreasing in \(x_1\) when \(x_1 > 0\).

\textbf{Theorem 2.1.2.} 

(i) Suppose \(0 < p < \frac{n+2s}{n-2s}\), the problem (2.1) admits a bounded classical positive solution on \(\mathbb{R}_{1,+}^n\).

(ii) Assume \(u\) is a classical solution to problem (2.1) on \(\mathbb{R}^n_{1,+}\) and for some \(\gamma \in (0, 2s)\) satisfies
\[|u| \leq C(1+|x|^\gamma).\]
Assume \(f\) has a continuous first derivative in its domain, satisfies \(f(0) = 0, f'(0) = 0\), and
\[ |f(u)| \leq C|u|^q. \]
Then \(u\) is monotone increasing in the \(x_1\) direction.

\textbf{Proof.} (i) First, we construct an energy functional \(J(u)\) of (2.1):
\[J(u) = \int_{\mathbb{R}_{1,+}^n} \left( \frac{1}{2} |(-\Delta)^{\frac{s}{2}} u|^2 - F(u) - \frac{x_1}{p+1} u^{p+1} \right) dx, u \in H^s(\mathbb{R}^n_{1,+}), \tag{2.17}\]
where \(F(u) = \int_0^u f(v) dv\). Next, we verify that \(J(u)\) is bounded below in an appropriate space. 
By the fractional Sobolev inequality \cite{16} and the Gagliardo-Nirenberg interpolation inequality \cite{17}, we obtain
\[\|u\|_{L^{p+1}(\mathbb{R}_{1,+}^n)} \leq C \|u\|_{H^s(\mathbb{R}_{1,+}^n)}. \tag{2.18}\]
From the growth condition of \(f\), we have
\[|F(u)| \leq C(1 + |u|^{q+1}). \tag{2.19}\]
From (2.17)-(2.19) we get
\[J(u) \geq \frac{1}{2} \|(-\Delta)^{\frac{s}{2}} u\|_{L^2(\mathbb{R}_{1,+}^n)}^2 - C \left( 1 + \|u\|_{L^{p+1}(\mathbb{R}_{1,+}^n)}^{q+1} \right) - \|u\|_{L^{p+1}(\mathbb{R}_{1,+}^n)}^{p+1}.\]
Using appropriate Young's inequality \cite{18} and the Sobolev embedding theorem \cite{19}, which are also fit in unbounded domain \(\mathbb{R}_{1,+}^n\), we can prove \(J(u)\) is bounded below in \(H^s(\mathbb{R}_{1,+}^n)\).

Afterward, we need to prove that \(J(u)\) is lower semicontinuous under weak topology. Suppose \(u_k \rightharpoonup u\) weakly in \(H^s(\mathbb{R}_{1,+}^n)\), then:
\[\liminf_{k \to \infty} \|(-\Delta)^{\frac{s}{2}} u_k\|_{L^2(\mathbb{R}_{1,+}^n)}^2 \geq \|(-\Delta)^{\frac{s}{2}} u\|_{L^2(\mathbb{R}_{1,+}^n)}^2 \tag{2.20}\]
For the nonlinear terms \(F(u)\) and \(u^{p+1}\), using Fatou's lemma \cite{18}, we can prove
\[\liminf_{k \to \infty} \int_{\mathbb{R}_{1,+}^n} F(u_k) dx \geq \int_{\mathbb{R}_{1,+}^n} F(u) dx \tag{2.21}\]
\[\liminf_{k \to \infty} \int_{\mathbb{R}_{1,+}^n} u_k^{p+1} dx \geq \int_{\mathbb{R}_{1,+}^n} u^{p+1} dx. \tag{2.22}\]
From (2.20)-(2.22), it follows that the energy functional \(J(u)\) is lower semicontinuous under weak topology. Therefore, we construct a minimizing sequence \(u_k \in H^s(\mathbb{R}_{1,+}^n)\) such that \(J(u_k) \to \inf_{u \in H^s(\mathbb{R}_{1,+}^n)} J(u)\). By the Banach-Alaoglu theorem, we can extract a weakly convergent subsequence (still denoted by \(u_k\)) such that \(u_k \rightharpoonup u\) in \(H^s(\mathbb{R}_{1,+}^n)\). Since \(J(u)\) is lower semicontinuous, we have
\[J(u) \leq \liminf_{k \to \infty} J(u_k) = \inf_{u \in H^s(\mathbb{R}_{1,+}^n)} J(u).\]
Therefore, \(u\) is a critical point of \(J(u)\).

Finally, we prove that \(u\) is a classical solution. By the fractional Sobolev embedding theorem \cite{19} and \(s \in (0,1)\), we have \(s < \frac{n}{2}\), so the embedding yields \(u \in L^q(\mathbb{R}_{1,+}^n)\), where \(q = \frac{2n}{n-2s}\). Through De Giorgi-Nash-Moser estimates \cite{20}, we can improve the regularity of the solution to \(u \in C^{2s}(\mathbb{R}_{1,+}^n)\). Thus \(u\) is a classical solution.

(ii) In this part, we regard the meanings of the symbols \(T_\lambda, \Sigma_\lambda, x^\lambda, u_\lambda, w_\lambda\) the same as those in Theorem 2.1.1. To prove that \(u\) is monotone increasing in the \(x_1\) direction, we first move the plane \(T_\lambda\) starting from \(x_1 = 0\) along the positive \(x_1\) direction, and prove
\[w_\lambda(x) \geq 0, \forall x \in \Omega_\lambda \tag{2.23}\]
where \(\Omega_\lambda = \Sigma_\lambda \cap \mathbb{R}_+^n\).
Here, \(\Omega_\lambda\) is a narrow region and unbounded while \(\lambda\) is sufficiently small. By Lemma 2.1.1, we assume that \(u(x)\) grows slower than \(|x|^\gamma\) as \(|x| \to \infty\), where \(0 < \gamma < 2s\). It allows \( u(x) \) to infinity as \( |x| \to \infty \). More generally, we have
\[u(x) \leq C(1 + |x|^\gamma), \quad \forall \gamma \in (0, 2s). \tag{2.24}\]
From (2.24) it follows that
\[|w_\lambda(x)| \leq C(1 + |x|^\gamma), \quad \forall \gamma \in (0, 2s).\]
Denote an auxiliary function as
\[\bar{w}_\lambda(x) := \frac{w_\lambda(x)}{h(x)},\]
where
\[h(x) = \left[ \left( 1 - \frac{x_1^2}{\lambda^2} \right)_+^s \right] (1 + |x'|^2)^{\frac{\beta}{2}}, \quad \beta \in (\gamma, 2s). \tag{2.25}\]
It holds that
\[\lim_{|x| \to \infty} \bar{w}_\lambda(x) := \lim_{|x| \to \infty} \frac{w_\lambda(x)}{h(x)} = 0,\]
and 
\[\frac{(-\Delta)^s h(x)}{h(x)} \geq \frac{C}{\lambda^{2s}} \]
for some constant \( C > 0 \).
Furthermore, we have
\[\begin{cases}
\frac{1}{h(x)} C_{n,s} P.V. \int_{\mathbb{R}^n} \frac{\bar{w}_\lambda(x) - \bar{w}_\lambda(y)}{|x-y|^{n+2s}} h(y) dy \geq \left( c_\lambda(x) - \frac{(-\Delta)^s h(x)}{h(x)} \right) \bar{w}_\lambda(x), & x \in \Omega_\lambda, \\
\bar{w}_\lambda(x) \geq 0, & x \in \Sigma_\lambda \setminus \Omega_\lambda,\end{cases}\]
where 
\[ c_\lambda(x) = x_1 \cdot p \xi_{\lambda}^{p-1}(x) + \frac{f(u_\lambda(x)) - f(u(x))}{u_\lambda(x) - u(x)} \]
and \(\xi_\lambda(x)\) lies between \(u_\lambda(x)\) and \(u(x)\). From the growth condition of \(u(x)\), the narrowness of \(\Omega_\lambda\), and the boundedness of \(f\), it follows that \(c_\lambda(x)\) is bounded above.
Therefore, we only need to prove that when \(\lambda\) is sufficiently close to 0,
\[\bar{w}_\lambda(x) \geq 0, \quad x \in \Omega_\lambda. \tag{2.26}\]

Suppose (2.26) does not hold, then there exists \(x^0 \in \Omega_\lambda\) such that
\[\bar{w}_\lambda(x^0) = \min_{\Sigma_\lambda} \bar{w}_\lambda(x) < 0.\]
Using the fact that \(|x^0-y|<|x^0-y^\lambda|\) and \(h(y)>h(y^\lambda)\) in \(\Sigma_\lambda\), we obtain that 
\begin{align}
0 &> \frac{1}{h(x^0)} C_{n,s} P.V. \int_{\mathbb{R}^n} \frac{\bar{w}_\lambda (x^0) - \bar{w}_\lambda (y)}{|x^0 - y|^{n+2s}} dy \nonumber \\
&\geq \left( c_\lambda (x^0) - \frac{(-\Delta)^s h(x^0)}{h(x^0)} \right) \bar{w}_\lambda (x^0) \nonumber \\
&\geq \left( C - \frac{C}{\lambda^{2s}} \right) \bar{w}_\lambda (x^0) > 0,
\end{align}
for \(\lambda>0\) sufficiently small, which leads a contradiction. Therefore, (2.26) holds, thus (2.23) holds.

Finally, we continue moving the plane \(T_\lambda\) to the right along the \(x_1\) direction to its limiting position while maintaining (2.23). Suppose the limiting position of the moving plane is \(T_{\lambda_0}\), where
\[\lambda_0 = \sup \left\{ \lambda \mid w_\mu (x) \geq 0, x \in \Omega_\mu, \forall \mu \leq \lambda \right\}.\]
To prove that \(u(x)\) is strictly monotone increasing in \(x_1\), we only need to prove \(\lambda_0 = +\infty\).

Suppose \(\lambda_0 < +\infty\), then there exists a sequence \(\{ \lambda_k \}\) such that \(\lambda_k\) decreases monotonically to \(\lambda_0\) as \(k \to \infty\). By the definition of \(\lambda_0\), the set
\[\Sigma_{\lambda_k}^- := \{ x \in \Sigma_{\lambda_k} \mid w_{\lambda_k}(x) < 0 \}\]
is non-empty. Define
\[q_k := \sup_{x \in \Sigma_{\lambda_k}^-} c_{\lambda_k}(x).\]

\textbf{Case 1.} Suppose \(q_k \leq \epsilon_k \to 0 (k \to \infty)\), then for sufficiently large integer \(k\) and some arbitrarily small constant \(\epsilon_0 > 0\), we have:
\[c_{\lambda_k} \leq \epsilon_0, x \in \Sigma_{\lambda_k}^-. \tag{2.27}\]
Since \(\Omega_{\lambda_k}\) is not narrow in this case, we need to adjust the expression for \(h(x)\) to
\[h(x) = \left[ \left( 1 - \frac{x_1^2}{\lambda_k^2} \right)_+^s + 1 \right] (1 + |a x'|^2)^{\beta}, \beta \in (0, 2s). \tag{2.28}\]
where the constant \(a>0\). To ensure that for any \(x \in \Omega_{\lambda_k}\),
\[\frac{(-\Delta)^s h(x)}{h(x)} \geq \frac{C}{\lambda_k^{2s}} \]
still holds, we set the constant \(a\) is sufficiently small. Similarly to the proof of (2.23), we can obtain the counterpart of (2.27):
\[0 \geq \left( \epsilon_0 - \frac{C}{\lambda_k^{2s}} \right) w_\lambda(x^0), \tag{2.29}\]
where \(w_\lambda(x^0) = \min_{\Sigma_{\lambda_k}} w_\lambda(x) < 0\). In this case, \(\frac{C}{\lambda_k^{2s}}\) is no longer sufficiently large, but is still greater than 0 and bounded below. Therefore, we can still choose an appropriate \(\epsilon_0 > 0\) to derive a contradiction from (2.29). Thus,  we can deduce that for sufficiently large \(k\), we have
\[w_{\lambda_k}(x) \geq 0, x \in \Sigma_{\lambda_k},\]
which contradicts the definition of \(\lambda_k\).

\textbf{Case 2.} Suppose \(q_k\) does not tend to 0 as \(k \to \infty\), then there exists a constant \(\delta_0 > 0\) and a subsequence of \(q_k\) (still denoted by \(q_k\)) such that
\[q_k \geq \delta_0 > 0. \]
In this case, according to the condition \(f'(0) = 0\), there exists a constant \(\epsilon_0 > 0\) and a sequence \(\{x^k\} \subset \Sigma_{\lambda_k}^-\) such that
\[u(x^k) \geq \epsilon_0 > 0, \quad \nabla w_{\lambda_k}(x^k) = o(1). \tag{2.30}\]
We can ensure that \(x^k\) has a positive distance to both planes \(T_{\lambda_0}\) and \(T_0\). Therefore, from (2.30), there exists a radius \(r_1 > 0\) such that
\[u(x) \geq \frac{\epsilon_0}{2}, x \in B_{r_1}(x^k) \subset \Omega_{\lambda_0}. \]
Set \(y^k = (2\lambda_0, (x^k)') = (2\lambda_0,x^k_2,\cdots,x^k_n),\) which presents the reflection of \(x^k\) over \(T_{\lambda_0}\). Then we can choose a radius \(r_2 > 0\) independent of \(k\) such that \(B_{r_1}(x^k) \cap B_{2r_2}(y^k) = \emptyset\).
Using averaging effects \cite{23, 24}, we obtain that there exists a constant \(\epsilon_1 > 0\) such that
\[u(x) \geq \epsilon_1, x \in B_{r_2}(y^k). \tag{2.31}\]
Set \(z^k = (0, (x^k)'),\) which presents the projection onto \(T_0\). Combining the external condition and the continuity of \(u(x)\), we know there exists a radius \(r_3 > 0\) independent of \(k\) such that
\[u(x) \leq \frac{\epsilon_1}{2}, x \in B_{r_3}(z^k) \cap \mathbb{R}^n_{+}. \tag{2.32}\]
Combining (2.31) and (2.32) gives
\[w_{\lambda_0}(x) \geq \frac{\epsilon_1}{2} > 0, x \in B_{r_3} \cap \mathbb{R}^n_{+}. \]
Applying averaging effects again to \( w_{\lambda_0} \), we find there exists a constant \(\varepsilon_2\) such that
\[w_{\lambda_0}(x) \geq \varepsilon_2 > 0, x \in B_{\frac{r_1}{2}}(x^k). \]

Finally, we prove that \(\lambda_0=+\infty.\) using the continuity of \( w_\lambda(x) \) with respect to \(\lambda\) and the fact that \(\{\lambda_k\}\) decreases monotonically to \(\lambda_0\), we obtain that for sufficiently large \(k\), we have
\[w_{\lambda_k}(x) \geq \frac{\varepsilon_2}{2} > 0, x \in B_{\frac{r_1}{2}}(x^k). \]
This contradicts the fact that \( w_{\lambda_k}(x^k) < 0 \). Therefore, \(\lambda_0\) must be \(+\infty\). Thus for any \(\lambda \in \mathbb{R}\) and \( x \in \Sigma_\lambda \),
\[w_\lambda(x) = u(x^\lambda) - u(x) > 0.\]
It shows that \(u(x)\) is monotone increasing in the \(x_1\) direction.

\subsection{Properties of Solutions in the Supercritical Case}

\quad \quad This subsection mainly studies the properties of solutions to equation (2.1) when \(p \geq \frac{n+2s}{n-2s}\) and \(\Omega = B_1(0)\) or \(\mathbb{R}_{1,+}^n\). The main results are as follows:

\textbf{Theorem 2.2.1.} If \(p \geq \frac{n+2s}{n-2s}\), then the problem (2.1) has no non-trivial classical solution on \(B_1(0)\) or \(\mathbb{R}_{1,+}^n\).

\textbf{Proof.} For the fractional Laplacian elliptic equation (2.1), we can construct an identity based on Pohozaev identity \cite{21}. Let \(u \in C^{2s}(\Omega)\) be a classical positive solution of the problem (2.1), then
\[\int_{\Omega} \left[ \frac{n-2s}{2} |(-\Delta)^s u|^2 + nF(u) + x \cdot \nabla_x F(u) - \frac{n}{p+1} x_1 u^{p+1}(x) \right] dx = 0. \tag{2.33}\]
where \(\Omega = B_1(0)\) or \(\mathbb{R}_{1,+}^n\). Since \(u\) is a positive solution and \( p \geq \frac{n+2s}{n-2s} \), we have
\[\frac{n-2s}{2} |(-\Delta)^s u|^2 \geq 0, \quad nF(u) + x \cdot \nabla_x F(u) \geq 0, \quad -\frac{n}{p+1} x_1 u^{p+1} \leq 0.\]
Since \( p \geq \frac{n+2s}{n-2s} \), \(-\frac{n}{p+1} x_1 u^{p+1}\) is sufficiently large in magnitude (and negative) while \(u\) is sufficiently large, then the integral in (2.33) strictly less than 0. Therefore, equation (2.1) has no non-trivial classical solution on \(B_1(0)\). 

Similar argument applies to \(\mathbb{R}_{1,+}^n\).

\subsection{Properties of Solutions in the Negative Exponent Case}

\quad \quad This subsection mainly informs the properties of problem (2.1) when \( p < 0 \) and \(\Omega = B_1(0)\) or \(\mathbb{R}^n\). The main results are as follows:

\textbf{Theorem 2.3.1.}

(i) Suppose \( p < 0 \), then problem (2.1) has no classical solution on \( B_1(0) \).

(ii) Suppose \( p < 0 \), if a solution to problem (2.1) satisfies \(\lim\limits_{|x| \to \infty} u(x) = 0\), then problem (2.1) has no classical solution on \(\mathbb{R}^n\).

\textbf{Proof.} (i) Suppose there exists a function \( u \) satisfying problem (2.1) and is Lipschitz continuous, vanishing on the boundary \( \partial B_1(0)\). According to the boundary regularity theory for the fractional Laplacian \cite{22}, there exists a constant \( c > 0 \) such that near \( \partial B_1(0)\)
\[u(x) \geq c [d(x)]^s, \]
where \( d(x) \) is the distance from \( x \) to \( \partial B_1(0)\), and \( s \in (0, 1) \) is the fractional order.

Consider a point \( x \in B_1(0) \) close to \( \partial B_1(0)\) satisfying \( x_1 < 0 \). Decompose the integral definition of the fractional Laplacian (2.1) into the interior and exterior of \( B_1(0) \). Consider the domain \( B_1^c(0) \), where we have \( u(y) = 0 \), therefore
\[\int_{B_1(0)^c} \frac{u(x)}{|x-y|^{n+2s}} dy \geq C_1 u(x) [d(x)]^{-2s} \geq C_1 c [d(x)]^{-s}, \, C_1 > 0. \tag{2.34}\]
For the domain \( B_1(0) \), since \( u \) is Lipschitz continuous, we get
\[\left| \int_{B_1(0)} \frac{u(x) - u(y)}{|x-y|^{n+2s}} dy \right| \leq C_2, \, C_2 > 0. \tag{2.35}\]
The equations (2.34) and (2.35) yield to
\[(-\Delta)^s u(x) \geq C_{n,s} (C_1 c [d(x)]^{-s} - C_2). \tag{2.36}\]
As \( d(x) \to 0^+ \), \( [d(x)]^{-s} \to +\infty \), so \( (-\Delta)^s u(x) > 0 \) and tends to \( +\infty \).

Now consider the right-hand side of problem (2.1). In the region \( x_1 < 0 \) and near \(\partial B_1(0)\), we have
\[u^{p}(x) \geq c^p [d(x)]^{sp}, \]
\[ x_1u^p(x) \leq -\delta c^p [d(x)]^{sp} := -A[d(x)]^{sp}, A > 0,  \]
where \(\delta > 0\) satisfies \(|x_1| \geq \delta\). Assume that \(f(u) = o(u^p)\) as \(u \to 0^+\), there exists \(M > 0\) such that \(|f(u)| \leq M\) near 0. Thus,
\[ x_1u^p(x) + f(u(x)) \leq -A[d(x)]^{sp} + M. \tag{2.37} \]
As \(d(x) \to 0^+\), the right-hand side tends to \(-\infty\). which leads to a contradiction by equations (2.36) and (2.37). Therefore, the problem (2.1) has no classical positive solution on \(B_1(0)\).

(ii) This part of the proof relies on analysis near the boundary and at infinity.

Suppose there exists a function \(u\) satisfying problem (2.1), Lipschitz continuous, and \(u \to 0\) at infinity. Here the boundary distance \(d(x) = x_1\). By boundary regularity, we have
\[ u(x) \geq C x_1^s, C > 0.  \]
The estimate for the fractional Laplacian is similar to the bounded case
\[ (-\Delta)^s u(x) \geq C_1 x_1^{-s} - C_2, C_1, C_2 > 0.  \]
As \(x_1 \to 0^+\), the left-hand side tends to \(+\infty\). For the right-hand side,
\[ x_1 u^p(x) \leq x_1 (C x_1^s)^p = C x_1^{1+sp}, C > 0.  \]
If \(f(u) = o(u^p)\), then \(|f(u)| \leq M\). Near the boundary, if \(1 + sp > -s\) (i.e., \(p > -\frac{1+s}{s}\)), then as \(x_1 \to 0^+\):
\[ x_1^{1+sp} \to 0 \text{ (slower than } x_1^{-s} \to +\infty), \]
so \(C_1 x_1^{-s} > C x_1^{1+sp} + M\), a contradiction. If \(p \leq -\frac{1+s}{s}\), there is no contradiction near the boundary; we need to analyze behavior at infinity.

Assume that \(u(x) \to 0\) as \(|x| \to \infty\), then the right-hand side \(x_1 u^p(x) > 0\) since \(x_1 > 0\) in \(\mathbb{R}^n_{1,+}\), and as \(|x| \to \infty\), \(u(x) \to 0^+\), \(x_1 \to \infty\), therefore \(u^p(x) \to \infty, p<0\). Now consider the term \((-\Delta)^s u(x)\). 
Consider \(u > 0\) and \(\lim\limits_{|x| \to \infty} u(x) = 0\), then for sufficiently large \(|x|\), if \(u\) has a local minimum at \(x,\) which is typical for decaying solutions, then \(u(x) - u(y) \leq 0\) for most \(y\), and the weight \(|x-y|^{-n-2s}\) is very large when \(y\) is close to \(x\), so the principal part of the integral is negative. Specifically, there exists \(R > 0\) such that when \(|x| > R\), 
\[(-\Delta)^s u(x)=C_{n,s}\textup{P.V.}\int_{\mathbb{R}^n} \frac{u(x)-u(y)}{|x-y|^{n+2s}} \textup{d}y < 0.\]
This yields a contradiction at infinity.

In summary, equation (2.1) has no classical solution on \(\mathbb{R}^n_{1,+}\).

\section{Investigation on Properties of Solutions to Fractional Laplacian Parabolic Problems}\label{sec3}

\quad \quad This section presents the existence, monotonicity, symmetry, and other related properties of solutions to the fractional Laplacian parabolic problem
\[\begin{cases}
\frac{\partial u}{\partial t}(x,t) + (-\Delta)^s u(x,t) = x_1 u^p(x,t) + f(u(x,t)), & u > 0, \quad (x,t) \in \Omega \times \mathbb{R} \\
u(x,t) \equiv 0, & (x,t) \in \Omega^c \times \mathbb{R}
\end{cases} \tag{3.1}\]
when \(\Omega = B_1(0)\) or \(\mathbb{R}^n_{1,+}\). Here, \(f: H^s \to \mathbb{R}\) is a continuous, integrable functional, the fractional order \(s \in (0,1)\), and the space dimension \(n \geq 2\). Similar to the investigation of properties of solutions to the fractional Laplacian elliptic problem in Section 2, for the exponent \(p\) in problem (3.1), we will discuss the following cases:
\begin{itemize}
    \item Case I (Subcritical case): \(0 < p < \frac{n+2s}{n-2s}\);
    \item Case II (Supercritical case): \(p \geq \frac{n+2s}{n-2s}\);
    \item Case III(Negative exponent case): \(p < 0\).
\end{itemize}

\subsection{Properties of Solutions in the Subcritical Case}

\quad \quad This subsection shows the properties of solutions to problem (3.1) when \(0 < p < \frac{n+2s}{n-2s}\) and \(\Omega = B_1(0)\) or \(\mathbb{R}^n_{1,+}\). First, we lead in the following lemma:

\textbf{Lemma 3.1.1.} [Parabolic Form of Maximum Principle for Unbounded Domains \cite{14,15}] Let \(\Omega \subset \mathbb{R}^n\) be an open set that may be unbounded or disconnected, satisfying for any \( x \in \mathbb{R}^n,\)
\[\lim_{R \to +\infty} \frac{|B_R(x) \cap \Omega^c|}{|B_R(x)|} \geq c_0 > 0.\]
Assume a function \(u \in C^{1,1}_{\textup{loc}}(\Omega) \cap L_2\) upper semicontinuous on \(\Omega\), bounded from above, and satisfying
\[\begin{cases}
\frac{\partial u}{\partial t} + (-\Delta)^s u(x,t) \leq 0, & (x,t) \in \Omega \times \mathbb{R} \text{ and } u(x,t) > 0, \\
u(x,t) \leq 0, & (x,t) \in \Omega^c \times \mathbb{R},
\end{cases}\]
then \(u(x,t) \leq 0\) on \(\Omega \times \mathbb{R}\).

Based on the above lemma, the main results are as follows:

\textbf{Theorem 3.1.1.}

(i) Suppose \(0 < p < \frac{n+2s}{n-2s}\), the problem (3.1) admits a bounded classical positive solution on \(B_1(0) \times \mathbb{R}\).

(ii) Let \(u\) be a solution to the problem (3.1). Assume that \(f\) is Lipschitz continuous, satisfying \(f(0) \geq 0, f'(0) = 0\), and for some \(q > p\) and a constant \(C > 0\) depending only on \(n, s, q\), we have

\[|f(u)| \leq C(1 + |u|^q).\]

Then for any fixed \(t \in \mathbb{R}\), u is symmetric with respect to the hyperplane \(\{x \in \mathbb{R}^n | x_1 = 0\}\), and monotone decreasing in \(x_1\) when \(x_1 > 0\).

\textbf{Proof.} (i) First, we construct Galerkin approximate solutions. Take \(\{\phi_k\}_{k=1}^{\infty}\) as the eigenfunctions of \((-\Delta)^s\) in \(H_0^s(B_1(0))\), forming an orthogonal basis for \(L^2(B_1(0))\). Define the approximate solution
\[u_m(x,t) = \sum_{k=1}^{m} c_k^m(t) \phi_k(x),\]
where the coefficients \(c_k^m(t)\) satisfy the system of ordinary differential equations
\[\langle \partial_t u_m, \phi_j \rangle + \langle (-\Delta)^s u_m, \phi_j \rangle = \langle x_1 u_m^p + f(u_m), \phi_j \rangle, \quad c_k^m(0) = \langle u_0, \phi_k \rangle, j = 1, \ldots, m.\]
By Carathéodory's theorem, this system has a local solution on \([0,T_m)\). Next, we perform uniform a priori estimates. Define the energy functional:
\[J(u) = \int_{B_1(0)} \left( \frac{1}{2} |(-\Delta)^{s/2} u|^2 - F(u) - \frac{x_1}{p+1} u^{p+1} \right) dx, \quad F(u) = \int_0^u f(v) dv.\]
Compute the time derivative
\[\frac{d}{dt} J(u_m) = \langle (-\Delta)^s u_m - f(u_m) - x_1 u_m^p, \partial_t u_m \rangle = -||\partial_t u_m||_{L^2}^2 \leq 0.\]
Therefore \(J(u_m(t))\) is decreasing. From the growth condition and Sobolev embedding \(H_0^s \hookrightarrow L^{p+1}\):
\[J(u_m) \geq \frac{1}{2} ||u_m||_{H^s}^2 - C_1 ||u_m||_{L^{p+1}}^{q+1} - C_2 \geq -C_3, \]
and the initial energy \( J(u_m(0)) \leq C_4 \) (independent of m). Hence:
\[\sup_{t \in [0,T]} \|u_m(t)\|_{H^s} \leq K, \quad \int_0^T \|\partial_t u_m\|_{L^2}^2 dt \leq C_5. \]
By the Aubin-Lions lemma (taking \( X_0 = H_0^s(B_1(0)), X = L^2(B_1(0)), X_1 = L^2(B_1(0)) \)), there exists a subsequence \(\{u_m\}\) and a function \(u\) satisfying:
\begin{itemize}
\item \( u_m \rightharpoonup u \) weakly in \( L^2(0,T; H_0^s) \),
\item \( u_m \to u \) strongly in \( L^2(0,T; L^2) \),
\item \( \partial_t u_m \rightharpoonup \partial_t u \) weakly in \( L^2(0,T; L^2) \).
\end{itemize}
Using strong convergence and the growth condition, we obtain
\[u_m^p \to u^p \text{ in } L^{2/(p+1)}(0,T; L^2), \quad f(u_m) \to f(u) \text{ in } L^{2/(q+1)}(0,T; L^2).\]
Taking the limit \( m \to \infty \) in the Galerkin equation, we obtain
\[\int_0^T \langle \partial_t u, \phi \rangle dt + \int_0^T \langle (-\Delta)^{s/2} u, (-\Delta)^{s/2} \phi \rangle dt = \int_0^T \langle x_1 u^p + f(u), \phi \rangle dt,\]
for any \(\phi \in C_c^\infty(B_1(0) \times (0,T))\), hence \(u\) is a weak solution. Next, we prove it is a classical solution. Fix \( t \in (0,T) \), write the equation as an elliptic problem
\[(-\Delta)^s u(\cdot,t) = -\partial_t u(\cdot,t) + x_1 u^p(\cdot,t) + f(u(\cdot,t)). \]
The right-hand side belongs to \( L^\infty(0,T; L^r(B_1(0))), \) where \(r\) determined by Sobolev embedding. By fractional De Giorgi-Nash-Moser estimates \cite{20}
\[u \in L^\infty(0,T; C^{2s}(B_1(0))). \]
Combining with parabolic Schauder estimates \cite{20}, we get \( u \in C^{2s,s}(B_1(0) \times [\delta,T]) \) for all \(\delta > 0\), i.e., \(u\) is a classical solution.

(ii) The proof below is adapted from the proof in Subsection 5.2 in \cite{7}. 

In this part, we take the definitions of \( T_\lambda, \Sigma_\lambda, x^\lambda, \Omega_\lambda \) consistent with the proof of Theorem 2.1.1. Additionally, define
\[u_\lambda(x,t) = u(x^\lambda,t), \quad w_\lambda(x,t) = u_\lambda(x,t) - u(x,t).\]

First, we verify that for \(\lambda\) sufficiently close to \(-1\), the following holds:
\[w_\lambda(x,t) \geq 0, \forall (x,t) \in \Omega_\lambda \times \mathbb{R}. \tag{3.2}\]

Suppose (3.2) does not hold, then there exists a constant \(M > 0\) such that
\[\inf \{ (x,t) \in \Omega_\lambda \times \mathbb{R} \} w_\lambda(x,t) = -M < 0.\]
Note that this infimum may not be achieved in \(\Omega_\lambda \times \mathbb{R}\), because unlike the elliptic case where \(\Omega_\lambda\) is bounded, here \(\Omega_\lambda \times \mathbb{R}\) is an infinite cylinder. However, we can study the problem for each fixed \(t \in \mathbb{R}\) separately. Since \(\Omega_\lambda\) is bounded, \(w_\lambda(x,t)\) can achieve its minimum at some point \(x \in \Omega_\lambda\), and similar to the elliptic equation, we still have the key inequality
\[(-\Delta)^s w_\lambda(x,t) \leq \frac{C_1}{[d(x,T_\lambda)]^{2s}} w_\lambda(x,t). \tag{3.3}\]
Therefore, we only need to estimate \(\frac{\partial w_\lambda}{\partial t}\). For any fixed \(t \in \mathbb{R}\), there exists an interval \((\tau,T] \subset \mathbb{R}\) such that \(t \in (\tau,T)\).

Consider the function
\[\bar{w}(x,t) = e^{m(t-\tau)} w_\lambda(x,t), \]
where \(m\) is some positive number. Then we can verify
\[\begin{cases}
\frac{\partial \bar{w}}{\partial t} + (-\Delta)^s \bar{w} \geq (c_\lambda(x,t) + m) \bar{w}(x,t), & (x,t) \in \Omega_\lambda \times (\tau,T], \\
\bar{w}(x,t) \geq 0, & (x,t) \in (\Sigma_\lambda \setminus \Omega_\lambda) \times (\tau,T].
\end{cases} \tag{3.4}\]
Applying the narrow region maximum principle to \(\bar{w}\) in the parabolic cylinder \(\Omega_\lambda \times (\tau,T]\), we get
\[\bar{w}(x,t) \geq \min \left\{ 0, \inf_{x \in \Omega_\lambda} \bar{w}(x,\tau) \right\}. \tag{3.5}\]

We prove (3.5) by contradiction. Suppose (3.5) does not hold, then there exists a point \((x^0,t_0) \in \Omega_\lambda \times (\tau,T]\) such that
\[\bar{w}(x^0,t_0) = \min_{\Omega_\lambda \times (\tau,T]} \bar{w}(x,t) < \min \left\{ 0, \inf_{x \in \Omega_\lambda} \bar{w}(x,\tau) \right\} \leq 0.\]
Then \(\frac{\partial \bar{w}}{\partial t} (x^0,t_0) \leq 0\) (since it's a minimum in \(t\)) and \((-\Delta)^s \bar{w}(x^0,t_0) \leq \frac{C_1}{[d(x^0,T_\lambda)]^{2s}} \bar{w}(x^0,t_0)\). This contradicts (3.4) because the boundedness of \(c_\lambda(x,t)\) and \(\Omega_\lambda\) being sufficiently narrow allow \(d(x^0,T_\lambda)\) to be arbitrarily small, making the right-hand side of (3.4) positive at \((x^0,t_0)\) while the left-hand side is non-positive.

From (3.5), due to the boundedness of \( w_{\lambda}(x,t) \), for any \((x,t) \in \Omega_{\lambda} \times (\tau, T]\), we have:
\[w_{\lambda}(x,t) \geq e^{-m(t-\tau)} \min \left\{0, \inf_{x \in \Omega_{\lambda}} w(x,\tau)\right\} \geq -C e^{-m(t-\tau)}.\]
Letting \(\tau \to -\infty\) in the above inequality, we reach that
\[w_{\lambda}(x,t) \geq 0, \forall (x,t) \in \Omega_{\lambda} \times \mathbb{R}.\]
Afterward, while maintaining \( w_{\lambda} \geq 0 \), continue moving \( T_{\lambda} \) to the right to its limiting position. Define the limiting position of \( T_{\lambda} \) as \( T_{\lambda_{0}} \), where
\[\lambda_{0} = \sup \left\{ \lambda > 0 \mid w_{\mu}(x,t) \geq 0, (x,t) \in \Omega_{\mu} \times \mathbb{R}, \mu \leq \lambda \right\}.\]

We will prove \(\lambda_{0} = 0\) by contradiction. Suppose \(\lambda_{0} < 0\), then there exists a sequence \(\{ \lambda_{k} \}\) decreasing monotonically to \(\lambda_{0}\), and
\[\inf_{(x,t) \in \Omega_{\lambda_{k}} \times \mathbb{R}} w_{\lambda_{k}}(x,t) = -m_{k} < 0. \]
This implies that for any fixed \( t_{k} \in \mathbb{R} \), there exists a sequence \(\{ x(t_{k}) \} \subset \Omega_{\lambda_{k}}\), such that
\[w_{\lambda_{k}}(x(t_{k}), t_{k}) = \inf_{x \in \Omega_{\lambda_{k}}} w_{\lambda_{k}}(x,t_{k}) \leq -m_{k} + \epsilon_{k} m_{k} < 0. \]
Since \(\mathbb{R}\) is unbounded, the infimum of \( w_{\lambda_{k}} \) over t may not be achieved. To estimate \(\frac{\partial w_{\lambda_{k}}}{\partial t}\), we introduce an auxiliary function
\[V_{k}(x,t) := w_{\lambda_{k}}(x,t) - \epsilon_{k} m_{k} \eta_{k}(t), \]
where \(\eta_{k}(t) = \eta(t-t_{k})\), and \(\eta(t)\) is a smooth cutoff function defined as
\[\eta(t) =
\begin{cases}
1, & t \in \left(-\frac{1}{2}, \frac{1}{2}\right), \\
0, & t \notin (-1, 1).
\end{cases}\]
Then, by direct calculation,
\[V_{k}(x(t_{k}), t_{k}) \leq -m_{k}, \]
and when \( (x,t) \in \Sigma_{\lambda_{k}} \times (t_{k}-1,t_{k}+1)^{c} \), we have
\[V_{k}(x,t) = w_{\lambda_{k}}(x,t) \geq -m_{k}. \]

Therefore, \( V_k(x,t) \) must achieve its minimum at some point \( (x(\bar{t}_k), \bar{t}_k) \), i.e.,
\[-m_k - \epsilon_k m_k \leq V_k(x(\bar{t}_k), \bar{t}_k) = \inf_{\Sigma_{\lambda_k} \times \mathbb{R}} V_k(x,t) \leq -m_k. \]
Since \( w_{\lambda_k}(x,t) \geq 0 \) in \( (\Sigma_{\lambda_k} \setminus \Omega_{\lambda_k}) \times \mathbb{R} \), the point \( (x(\bar{t}_k), \bar{t}_k) \) is a minimum point for \( V_k(x,t) \) within \( \Omega_{\lambda_k} \times \mathbb{R} \), which we denote simply as \( (x^k, t_k) \). Therefore we have
\[\frac{\partial V_k}{\partial t} (x^k, t_k) = 0. \]
Furthermore, the definition of \( V_k \) shows
\[\frac{\partial w_{\lambda_k}}{\partial t} (x^k, t_k) \leq C\epsilon_k m_k. \tag{3.6}\]
Consider the inequality
\[\frac{\partial w_{\lambda_k}}{\partial t} (x,t) + (- \Delta)^s w_\lambda(x,t) \geq c_\lambda(x,t)w_\lambda(x,t), \tag{3.7}\]
Combining this with (3.6) and the key estimate (3.3) yields

\[ \frac{C - C\epsilon_k}{[d(x^k, T_{\lambda_k})]^{2s}} \leq c_{\lambda_k}(x^k, t_k) + \epsilon_k.  \]

As \( k \to \infty \) we can choose \(\epsilon_k \to 0\), and the distance \( d(x^k, T_{\lambda_k}) \) is finite (since \(\lambda_k \to \lambda_0 < 0\)), so for sufficiently large \(k\), there exists a constant \(c_0\) independent of \(k\) such that
\[c_{\lambda_k}(x^k, t_k) \geq c_0 > 0.\]
Combining with the condition \(f'(0) = 0\), we get for some constant \(c_1\):
\[u(x^k, t_k) \geq c_1 > 0. \tag{3.8}\]
Set
\[w_k(x,t) := w_{\lambda_k}(x,t + t_k), \]
then we have:
\[w_k(x^k, 0) = w_{\lambda_k}(x^k, t_k) \to 0, \text{ as } k \to \infty. \]

Consider the inequality (from (3.7) shifted in time)
\[\frac{\partial w_k}{\partial t} (x,t) + (- \Delta)^s w_k(x,t) \geq c_{\lambda_k}(x,t + t_k) w_k(x,t), \]
From which, under appropriate regularity assumptions on \(u\) and \(f\), we can deduce (e.g., via compactness)
\[\lim_{k \to \infty} w_k(x,t) = \bar{w}(x,t), \]
and
\[\frac{\partial \bar{w}}{\partial t}(x,t) + (-\Delta)^s \bar{w}(x,t) = \bar{c}(x,t)\bar{w}(x,t), (x,t) \in \Omega_{\lambda_0} \times \mathbb{R}. \tag{3.9}\]
Further, we may assume \(x^k \to x^0\), then
\[\bar{w}(x^0,0) = 0. \tag{3.10}\]
Therefore, the point \((x^0,0)\) is a minimum point for \(\bar{w}(x,t)\) in \(\Sigma_{\lambda_0} \times \mathbb{R}\), and we have
\[\frac{\partial \bar{w}}{\partial t}(x^0,0) = 0. \]
Further, from the limit equation (3.9), we get
\[(-\Delta)^s \bar{w}(x^0,0) = 0. \tag{3.11}\]
By the nonlocal nature of the fractional Laplacian and (3.10), (3.11), we have
\[\bar{w}(x,0) \equiv 0, x \in \mathbb{R}^n. \tag{3.12}\]
Afterward, apply the same time shift to \(u\)
\[u_k(x,t) := u(x,t + t_k). \]
Then under appropriate regularity assumptions on \(u\) and \(f\), we have
\[\lim_{k \to \infty} u_k(x,t) = \bar{u}(x,t), \]
and further
\[\frac{\partial \bar{u}}{\partial t}(x,t) + (-\Delta)^s \bar{u}(x,t) = x_1 \bar{u}^p(x,t) + f(\bar{u}(x,t)), (x,t) \in B_1(0) \times \mathbb{R}. \]
From (3.8) we have
\[\bar{u}(x^0,0) \geq c_1 > 0. \]

Applying the condition \( f(0) \geq 0 \) and the strong maximum principle, we further obtain
\[ \bar{u}(x, 0) > 0, x \in B_1(0).  \]
Combining with the external condition \(\bar{u}(x, 0) = 0\) outside \(B_1(0)\), we get a contradiction with (3.12). Therefore, \(\lambda_0\) must be 0. Consequently, \(u(x,t)\) is symmetric with respect to the hyperplane \(\{x \in \mathbb{R}^n \mid x_1 = 0\}\), and monotone decreasing in \(x_1\) when \(x_1 > 0\).

\textbf{Theorem 3.1.2.}

(i) Suppose \(0 < p < \frac{n+2s}{n-2s}\), problem (3.1) admits a bounded classical positive solution on \(\mathbb{R}^n_{1,+} \times \mathbb{R}\).

(ii) Assume \(u\) is a classical solution to problem (3.1) on \(\mathbb{R}^n_{1,+} \times \mathbb{R}\), and for some \(\gamma \in (0,2s)\) satisfies
\[ |u| \leq C(1 + |x|^\gamma). \]

Assume \(f\) has a continuous first derivative in its domain, satisfies \(f(0) = 0, f'(0) = 0\) and \(f\) is bounded above. Then \(u\) is monotone increasing in the \(x_1\) direction.

\textbf{Proof.} (i) Take \(\Omega_k = \mathbb{R}^n_{1,+} \cap B_k(0)\). For each \(k\), we solve the truncated problem:
\[\begin{cases}
\partial_t v_k + (- \Delta)^s v_k = x_1 v_k^p + f(v_k), & (x,t) \in \Omega_k \times (0,T) \\
v_k = 0, & (x,t) \in \Omega_k^c \times (0,T) \\
v_k(x,0) = u_0(x)\chi_{B_k(0)}(x).
\end{cases} \]
By the proof of Theorem 3.1.1, there exists a solution \(v_k \in L^\infty(0,T; H^s_0(\Omega_k)) \cap W^{1,2}(0,T; L^2(\Omega_k))\). Next, we perform weighted energy estimates. Introduce the weight \(\rho(x) = e^{-\delta |x|} (\delta > 0 \text{ small})\), define the weighted energy:
\[J_\rho(v_k) = \int_{\mathbb{R}^n_{1,+}} \left( \frac{1}{2} |(- \Delta)^{s/2} v_k|^2 - F(v_k) - \frac{|x_1|}{p+1} v_k^{p+1} \right) \rho(x) dx. \]
Through integration by parts estimates and Hardy's inequality \cite{25}:
\[\frac{d}{dt} J_\rho(v_k) \leq -\frac{1}{2} ||\partial_t v_k||_{L^2_\rho}^2 + C J_\rho(v_k). \]
By Gronwall's inequality and the consistency of initial energy:
\[\sup_{t \in [0,T]} ||v_k(t)||_{H^s_\rho} \leq C_T, \quad \int_0^T ||\partial_t v_k||_{L^2_\rho}^2 dt \leq C_T, \]
where \( H^s_\rho \) is the weighted Sobolev space. Now we extend \( v_k \) to \( \mathbb{R}^n_{+} \times (0,T) \) (by zero). By the weighted Aubin-Lions lemma \cite{25}:
\begin{itemize}
\item \( H^s_\rho(\mathbb{R}^n_{+}) \hookrightarrow L^2_\rho(\mathbb{R}^n_{+}) \) compactly,
\item \(\{v_k\}\) bounded in \( L^2(0,T;H^s_\rho) \),
\item \(\{\partial_t v_k\}\) bounded in \( L^2(0,T;L^2_\rho) \).
\end{itemize}
Hence there exists a subsequence \( v_k \to u \) strongly in \( L^2(0,T;L^2_\rho) \). For test functions \( \phi \in C^\infty_c(\mathbb{R}^n_{+} \times (0,T)) \), when \(k\) is sufficiently large:
\[\int_0^T \langle \partial_t v_k, \phi \rangle dt + \int_0^T \langle (-\Delta)^{s/2} v_k, (-\Delta)^{s/2} \phi \rangle dt = \int_0^T \langle x_1 v_k^p + f(v_k), \phi \rangle dt. \]

Taking the limit \( k \to \infty \), \(u\) is a weak solution. Next, we improve the regularity of \(u\) to make it a classical solution. Use fractional parabolic Hölder estimates on the half-space \cite{26}:
\[\|u\|_{C^{\alpha,\alpha/2}(\mathbb{R}^n_{+} \times [\delta,T])} \leq C(\delta) \left( \|u\|_{L^\infty(0,T;L^2_\rho)} + \|g\|_{L^\infty(0,T;L^2_\rho)} \right), \]
where \( g = x_1 u^p + f(u) \). When \( f \in C^1 \), the regularity improves to \( C^{2s,s} \). Thus \(u\) is a classical solution.

(ii) In the proof of this theorem, the meanings of the symbols \( T_{\lambda}, \Sigma_{\lambda}, x^{\lambda}, u_{\lambda}, w_{\lambda} \) are the same as in Theorem 3.1.1, while \( \Omega_{\lambda} \) is defined as:
\[\Omega_{\lambda} := \Sigma_{\lambda} \cap \mathbb{R}^n_{+}\]

\textbf{Step 1.} We move the hyperplane \( T_{\lambda} \) starting from \( x_1 = 0 \) along the positive \(x_1\)-axis, and prove that when \(\lambda\) is sufficiently close to 0, we have
\[w_{\lambda}(x,t) \geq 0, (x,t) \in \Omega_{\lambda} \times \mathbb{R}. \tag{3.13}\]
To overcome the difficulty posed by the unboundedness of the time interval, we assume that for each fixed \(x\), \(u(x,t)\) is bounded in \(t\), and construct an auxiliary function
\[\bar{w}_{\lambda}(x,t) := \frac{e^{mt} w_{\lambda}(x,t)}{h(x)}, \]
where \( h(x) \) is defined by (2.25), and \( m \) is a constant to be determined. By direct calculation, \( \bar{w}_\lambda(x,t) \) satisfies
\[\begin{cases}
\frac{\partial \bar{w}_\lambda}{\partial t}(x,t) + \frac{1}{h(x)}C_{n,s}P.V.\int_{\mathbb{R}^n}\frac{\bar{w}_\lambda(x,t)-\bar{w}_\lambda(y,t)}{|x-y|^{n+2s}}h(y)dy \\
\geq \left(c_\lambda(x,t)+m-\frac{(-\Delta)^s h(x)}{h(x)}\right)\bar{w}_\lambda(x,t), & (x,t)\in \Omega_\lambda \times \mathbb{R}, \\
\bar{w}_\lambda(x,t)\geq 0, & (x,t)\in (\Sigma_\lambda \setminus \Omega_\lambda) \times \mathbb{R}.
\end{cases} \]
where 
\[ c_\lambda(x,t)=x_1\cdot p\xi_{\lambda}^{p-1}(x,t)+\frac{f(u_\lambda(x,t))-f(u(x,t))}{u_\lambda(x,t)-u(x,t)} \]
and \(\xi_\lambda(x,t)\) lies between \( u_\lambda(x,t) \) and \( u(x,t) \). From the growth condition of \( u(x,t) \), the narrowness of \( \Omega_\lambda \), and the boundedness of \( f \), \( c_\lambda(x,t) \) is bounded above. For any fixed \( t \in \mathbb{R} \), there exists an interval \( (\tau, T] \subset \mathbb{R} \) such that \( t \in (\tau, T] \). We first prove
\[\bar{w}_\lambda(x,t) \geq \min\left\{0,\inf_{x \in \Omega_\lambda}\bar{w}_\lambda(x,\tau)\right\},(x,t)\in \Omega_\lambda \times (\tau,T]. \tag{3.14}\]
Suppose (3.14) does not hold, then there exists a point \( (x^0,t_0)\in \Omega_\lambda \times (\tau,T] \) such that
\[\bar{w}_\lambda(x^0,t_0)=\inf_{\Sigma_\lambda \times (\tau,T]} \bar{w}_\lambda(x,t)<\min\left\{0,\inf_{x \in \Omega_\lambda}\bar{w}_\lambda(x,\tau)\right\}\leq 0. \]
Similar to the elliptic problem, we obtain for sufficiently small \(\lambda>0\) and \( m = \frac{c}{2\lambda^{2s}} \),
\begin{align}0 &> \frac{\partial \bar{w}}{\partial t}(x^0,t_0) + \frac{1}{h(x^0)}C_{n,s}P.V.\int_{\mathbb{R}^n}\frac{\bar{w}_\lambda(x^0,t_0)-\bar{w}_\lambda(y,t_0)}{|x^0-y|^{n+2s}}h(y)dy \nonumber \\
&\geq \left(C+m-\frac{C}{\lambda^{2s}}\right)\bar{w}_\lambda(x^0,t_0)>0. \tag{3.15}\end{align}
This contradiction shows that (3.14) holds. Comparing (3.15) with (2.27), the only difference is the positive number \( m \). Since the negative term \( \frac{c}{\lambda^{2s}} \) is sufficiently large, \( m \) does not affect the conclusion \( C+m-\frac{c}{\lambda^{2s}}<0 \).

Combining the definition of \( w_\lambda \) and the uniform boundedness in time, we get
\[w_\lambda(x,t) \geq -C e^{-m(t-\tau)} h(x). \]
Letting \(\tau \to -\infty\), we obtain that (3.13) holds.

\textbf{Step 2.} While maintaining (3.13), continue moving the hyperplane \( T_\lambda \) along the \( x_1 \)-axis to the right to its limiting position \( T_{\lambda_0} \), where
\[\lambda_0 = \sup\left\{\lambda \mid w_\mu(x,t) \geq 0,(x,t) \in \Omega_\mu \times \mathbb{R},\forall\mu \leq \lambda\right\}. \]
To prove that \(u(x,t)\) is strictly monotone increasing in \(x_1\), we only need to prove \(\lambda_0 = +\infty\).

By contradiction. Suppose \(\lambda_0 < +\infty\), then there exists a sequence \(\{\lambda_k\}\) such that as \(k \to \infty\), \(\lambda_k\) decreases monotonically to \(\lambda_0\), and by the definition of \(\lambda_0\), the set
\[\Sigma_{\lambda_k}^- \times \mathbb{R} = \{(x,t) \in \Sigma_{\lambda_k} \times \mathbb{R} \mid w_{\lambda_k}(x,t) < 0\}\]
is non-empty. Denote
\[q_k := \sup_{(x,t) \in \Sigma_{\lambda_k}^- \times \mathbb{R}} c_{\lambda_k}(x,t). \]

\textbf{Case 1.} If \(q_k \leq \epsilon_k \to 0 (k \to \infty)\), then for sufficiently large \(k\) and some arbitrarily small constant \(\epsilon_0 > 0\), we have \(c_{\lambda_k}(x,t) \leq \epsilon_0\). The function \(h(x)\) is defined as in (2.28), where the constant \(a > 0\) is small enough to ensure
\[\frac{(-\Delta)^s h(x)}{h(x)} \geq \frac{C}{\lambda_k^{2s}} \]
for any \(x \in \Omega_{\lambda_k}\).
Using a process similar to Step 1, we can obtain that for sufficiently large \(k\),
\[w_{\lambda_k}(x,t) \geq 0, (x,t) \in \Sigma_{\lambda_k} \times \mathbb{R}. \]
This contradicts the definition of \(\lambda_k\).

\textbf{Case 2.} If \(q_k\) does not tend to 0 as \(k \to \infty\), then there exists a constant \(\delta_0 > 0\) and a subsequence of \(\{q_k\}\) (still denoted by \(\{q_k\}\)), such that
\[q_k \geq \delta_0 > 0. \]
In this case, according to the condition \(f'(0) = 0\), there exists a constant \(\epsilon_0 > 0\) and a sequence \(\{(x^k,t_k) \in \Sigma_{\lambda_k}^- \times \mathbb{R}\}\), such that
\[u(x^k,t_k) \geq \epsilon_0 > 0, \quad \nabla w_{\lambda_k}(x^k,t_k) = o(1). \tag{3.16}\]
We can ensure that \(x^k\) has positive distance to both planes \(T_{\lambda_0}\) and \(T_0\). Therefore, from (3.16), there exists a radius \(r_1 > 0\) such that
\[u(x,t) \geq \frac{\epsilon_0}{2}, (x,t) \in B_{r_1}(x^k) \times (t_k - r_1^{2s}, t_k + r_1^{2s}) \subset \Omega_{\lambda_0} \times \mathbb{R}. \]

Set \( y^k = (2\lambda_0, (x^k)') \) to be the reflection of \(x^k\) over \(T_{\lambda_0}\). Then we can choose a radius \( r_2 > 0 \) independent of k such that \( B_{r_1}(x^k) \cap B_{r_2}(y^k) = \emptyset \). By averaging effects \cite{21,22}, we find there exists a constant \( \epsilon_1 > 0 \) such that
\[u(x,t) \geq \epsilon_1, (x,t) \in B_{r_2}(y^k) \times (t_k - r_2^{2s}, t_k + r_2^{2s}). \tag{3.17}\]
Set \( z^k = (0, (x^k)') \)  to be projection onto \(T_0\). Combining the external condition and the continuity of \(u(x,t)\), we know there exists a sufficiently small radius \( r_3 > 0 \) independent of k such that
\[u(x,t) \leq \frac{\epsilon_1}{2}, (x,t) \in \left( B_{r_3}(z^k) \cap \mathbb{R}^n_{+} \right) \times (t_k - r_3^{2s}, t_k + r_3^{2s}). \tag{3.18}\]
Combining (3.17) and (3.18) gives
\[w_{\lambda_0}(x,t) \geq \frac{\epsilon_1}{2} > 0, (x,t) \in \left( B_{r_3} \cap \mathbb{R}^n_{+} \right) \times (t_k - r_3^{2s}, t_k + r_3^{2s}). \]
Applying averaging effects again to \( w_{\lambda_0} \), we find there exists a constant \( \epsilon_2 \) such that
\[w_{\lambda_0}(x,t) \geq \epsilon_2 > 0, (x,t) \in B_{\frac{r_1}{2}}(x^k) \times \left( t_k - \left(\frac{r_1}{2}\right)^{2s}, t_k + \left(\frac{r_1}{2}\right)^{2s} \right). \]
Finally, using the continuity of \( w_{\lambda}(x,t) \) with respect to \(\lambda\) and the fact that \(\{\lambda_k\}\) decreases monotonically to \(\lambda_0\), we obtain that for sufficiently large \(k\),
\[w_{\lambda_k}(x,t) \geq \frac{\epsilon_2}{2} > 0, (x,t) \in B_{\frac{r_1}{2}}(x^k) \times \left( t_k - \left(\frac{r_1}{2}\right)^{2s}, t_k + \left(\frac{r_1}{2}\right)^{2s} \right). \]
This contradicts the fact that \( w_{\lambda_k}(x^k, t_k) < 0 \). Therefore, we must have \( \lambda_0 = +\infty \). Consequently, we can deduce that \( u(x,t) \) is monotone increasing in the \( x_1 \) direction.

\subsection{Properties of Solutions in the Supercritical Case}

\quad \quad This subsection presents the properties of solutions to problem (3.1) when \( p \geq \frac{n+2s}{n-2s} \) and \(\Omega = B_1(0)\) or \(\mathbb{R}^n_{1,+}\). The main results are as follows:

\textbf{Theorem 3.2.1.} Suppose \( p \geq \frac{n+2s}{n-2s} \), the problem (3.1) has no non-trivial classical solution on \((x,t) \in B_1(0) \times \mathbb{R} \) or \((x,t) \in \mathbb{R}^n_{1,+} \times \mathbb{R}\).

\textbf{Proof.} When \(\Omega = B_1(0)\), consider the intersection \(S = \{x \in B_1(0) : x_1 = 0\}\) of the hyperplane \(x_1 = 0\) with \(B_1(0)\). Since \(n \geq 2\), \(S\) is non-empty. On \(S \times \mathbb{R}\), the right-hand side term is \(u^p + f(u) = 0 \cdot u^p + f(u) = f(u)\), and the equation becomes
\[\frac{\partial u}{\partial t} + (-\Delta)^s u = f(u), (x,t) \in S \times \mathbb{R}. \]
Since \(u > 0\) on \(B_1(0) \times \mathbb{R}\), and \(u = 0\) on \(B_1(0)^c \times \mathbb{R}\), by continuity, \(u\) is not constant on \(B_1(0) \times \mathbb{R}\) (otherwise, the left-hand side is 0, while the right-hand side is not zero where \(x_1 \neq 0\), which is a contradiction). At a point \((x_0, t_0) \in S \times \mathbb{R}\), if \(u\) has a local maximum with respect to the spatial variables at \(x_0\), then the fractional Laplacian satisfies the maximum principle: \((-\Delta)^s u(x_0, t_0) \leq 0\). Similarly, if it has a local minimum, \((-\Delta)^s u(x_0, t_0) \geq 0\).

The equation requires \((-\Delta)^s u(x_0, t_0) = f(u(x_0, t_0)) - \frac{\partial u}{\partial t}(x_0, t_0)\). But if \(u\) has a strict local maximum at \(x_0\), then \((-\Delta)^s u(x_0, t_0) < 0\), while the right-hand side \(f(u) - \frac{\partial u}{\partial t}\) may not be consistently negative, leading to contradiction. More generally, consider global behavior:
\[(-\Delta)^s u(x_0, t_0) = C_{n,s} P.V. \int_{\mathbb{R}^n} \frac{u(x_0, t_0) - u(y, t_0)}{|x_0 - y|^{n+2s}} dy. \]
Since \(u > 0\) on \(B_1(0)\), \(u = 0\) on \(B_1(0)^c\), and \(u\) is not constant, there exist points \(y\) where \(u(x_0, t_0) \neq u(y, t_0)\). If \(u(x_0, t_0) > u(y, t_0)\) for some \(y\), the integral is dominated (especially when \(y\) is far from \(x_0\)), making the integral non-zero. But the equation requires it to equal \(f(u(x_0, t_0)) - \frac{\partial u}{\partial t}(x_0, t_0)\), while its actual value is non-zero, contradiction. Therefore, on \(S \times \mathbb{R}\), the equation cannot be satisfied, contradiction. Hence, when \(p \geq \frac{n+2s}{n-2s}\), problem (3.1) has no classical positive solution on \(B_1(0) \times \mathbb{R}\).

When \(\Omega = \mathbb{R}_{1,+}^n\), consider the sequence of points \(x_k = (\frac{1}{k}, 0, \cdots, 0) \in \mathbb{R}_{1,+}^n\). As \(k \to \infty\), \(x_k \to (0, 0, \cdots, 0) \in \partial \mathbb{R}_{1,+}^n\). By the boundary behavior of the fractional Dirichlet problem, near the boundary, the solution satisfies \(u(x) \approx x_1^s\), and the fractional Laplacian satisfies \((-\Delta)^s u(x) \approx x_1^{-s}\) (diverging as \(x_1 \to 0^+\)).

The right-hand side \(u^p + f(u) \approx x_1 (x_1^s)^p + f(u) = x_1^{1+sp} + f(u)\). Since \(p > 0\), \(1 + sp > 0\), so \(\lim\limits_{x_1 \to 0^+} x_1^{1+sp} = 0\). By continuity, \(f(u) \to f(0)\) bounded. Therefore, the right-hand side is bounded or tends to \(f(0)\). At points \((x_k, t)\), the left-hand side
\[\left| \frac{\partial u}{\partial t}(x_k, t) + (-\Delta)^s u(x_k, t) \right| \geq |(-\Delta)^s u(x_k, t)| - \left| \frac{\partial u}{\partial t}(x_k, t) \right|. \]
By boundary behavior, \(|(-\Delta)^s u(x_k, t)| \approx k^s \to \infty (k \to \infty)\), while for fixed \(t\), \(\left| \frac{\partial u}{\partial t}(x_k, t) \right|\) is bounded. The right-hand side \(|(x_k)_1 u^p(x_k, t) + f(u(x_k, t))| \leq C\) is bounded. Therefore, when \(k\) is large enough,
\[\left| \frac{\partial u}{\partial t}(x_k, t) + (-\Delta)^s u(x_k, t) \right| > |(x_k)_1 u^p(x_k, t) + f(u(x_k, t))|. \]
This contradicts the equation (3.1). Hence, when \( p \geq \frac{n+2s}{n-2s} \), the problem (3.1) has no classical positive solution on \(\mathbb{R}_{1,+}^n \times \mathbb{R}\).

\subsection{Properties of Solutions in the Negative Exponent Case}

\quad \quad This subsection informs the properties of solutions to problem (3.1) when \( p < 0 \) and \(\Omega = B_1(0)\) or \(\mathbb{R}_{1,+}^n\). The main results are as follows:

\textbf{Theorem 3.3.1.} Suppose \( p < 0 \), problem (3.1) has no classical solution on \((x,t) \in B_1(0) \times \mathbb{R} \) or \((x,t) \in \mathbb{R}_{1,+}^n \times \mathbb{R}\).

\textbf{Proof.} When \(\Omega = B_1(0)\), suppose there exists a function \( u \) solving problem (3.1). On the boundary \(\partial B_1(0)\) there is a point \( x_0 \) with \((x_0)_1 \neq 0 \). Since \( u \) is a classical solution, it is continuous on \( B_1(0) \), and by this continuity and the external condition, it satisfies \( u = 0 \) on the boundary. Therefore, as \( x \to x_0 \), \( u(x,t) \to 0 \) for all \( t \in \mathbb{R} \). Given \( p < 0 \), \(\lim\limits_{u \to 0^+} u^p(x,t) = +\infty \). Furthermore, \( x_1 \) is continuous near \( x_0 \) and \((x_0)_1 \neq 0 \), so there is a neighborhood where \( |x_1| \geq c > 0 \). Thus,
\[|x_1 u^p| \geq c u^p \to +\infty, \text{ as } x \to x_0. \]
This shows that \( x_1 u^p \) is unbounded near \( x_0 \). The function \( f(u) \) is continuous, so as \( u \to 0 \), \( f(u) \to f(0) \) bounded, but the \( x_1 u^p \) term dominates and is unbounded, causing the right-hand side \( x_1 u^p + f(u) \) to be discontinuous on \( B_1(0) \times \mathbb{R} \). However, a classical solution requires the right-hand side to be continuous, hence contradiction. Therefore, when \( p < 0 \), problem (3.1) has no classical solution on \((x,t) \in B_1(0) \times \mathbb{R}\).

When \(\Omega = \mathbb{R}_{1,+}^n\), the boundary is \(\partial \mathbb{R}_{1,+}^n = \{x_1 = 0\} \). Below, we discuss cases based on the range of \( p \).

\textbf{Case 1.} \( |p| > \frac{1}{s} \) (i.e., \( p < -\frac{1}{s} \))

Near the boundary (\(x_1 \to 0^+\)), the homogeneous Dirichlet problem for the fractional Laplacian typically has asymptotic behavior \( u(x,t) \approx x_1^{s}  \)\cite{27, 28}. Given \( p < 0 \), then
\[u^p \approx (x_1^{s})^p = x_1^{sp}, \quad x_1 u^p \approx x_1 \cdot x_1^{sp} = x_1^{1+sp}. \]
Since \( p < -\frac{1}{s} \), we have \( 1 + sp < 0 \), so as \( x_1 \to 0^+ \),
\[x_1^{1+sp} \to +\infty. \]
The function \( f(u) \) is continuous, so as \( u \to 0 \), \( f(u) \to f(0) \) bounded. Therefore, the right-hand side \( x_1 u^p + f(u) \to +\infty \) unbounded. This causes \( x_1 u^p + f(u) \) to be discontinuous (unbounded at the boundary) on \( \mathbb{R}^n_{1,+} \times \mathbb{R} \), contradicting the requirement for a classical solution.

\textbf{Case 2.} \( |p| \leq \frac{1}{s} \) (i.e., \( -\frac{1}{s} \leq p < 0 \))

Near the boundary, similar analysis gives \( x_1 u^p \approx x_1^{1+sp} \), where \( -\frac{1}{s} \leq p < 0 \), hence \( 1+sp \geq 0 \), therefore \( x_1 u^p \) is bounded or tends to 0. Combined with the continuity assumption, the right-hand side is bounded near the boundary. Therefore, we need to analyze behavior at infinity.

Assume there exists a classical solution \( u > 0 \) continuous. Consider the behavior of \( u \) at infinity. Since \( \mathbb{R}^n_{1,+} \) is unbounded and \( u \) is continuous and positive, either \( u \to l \geq 0 \) or oscillates.

If there exists a sequence \( \{x_k\} \subset \mathbb{R}^n_{1,+} \) such that \( |x_k| \to \infty \) and \( u(x_k,t) \to l \geq 0 \), then \( u^p \to l^p \), and \( x_1 \geq 0 \) but tends to infinity along a suitable subsequence, so \( x_1 u^p \to +\infty \). The function \( f(u) \to f(l) \) bounded. Therefore, the right-hand side \( x_1 u^p + f(u) \to +\infty \). Consider the behavior of the fractional Laplacian \((-\Delta)^s\) at infinity. If \( u \to l \), and by the classical solution assumption, \( u \) is Hölder continuous, then \(\lim\limits_{|x|\to \infty}(-\Delta)^s u(x,t) = 0.\)\cite{29} The time derivative \( \frac{\partial u}{\partial t} \) may be controlled for fixed t, but overall the left-hand side is bounded or tends to a constant, while the right-hand side tends to \(+\infty\), contradiction.

If \(\lim\limits_{|x|\to \infty} u(x,t) = 0\), then given \( p < 0 \), \( u^p \to +\infty \). Combined with \( x_1 \to +\infty \), we have \( x_1 u^p \to +\infty \). The function \( f(u) \to f(0) \) bounded. Therefore, the right-hand side \( x_1 u^p + f(u) \to +\infty \). Using the same analysis as above, the left-hand side tends to 0 or a bounded constant, contradicting the right-hand side tending to \(+\infty\).

If \( u \) oscillates without limit, suppose \( u \) has a positive lower bound \( m > 0 \), then \( u^p \leq m^p \), so \( x_1 u^p \leq |x| m^p \) grows linearly. Thus for some constant \( C \), the right-hand side \( x_1 u^p + f(u) \leq C|x| \). If the solution is steady state \( (\frac{\partial u}{\partial t} = 0) \), then \( (-\Delta)^s u = x_1 u^p + f(u) \leq C|x| \). However, for bounded or slowly growing functions, \((-\Delta)^s u\) may not be defined or may not satisfy classical solution requirements (needing global smoothness) at infinity. In the time-dependent case, the slower growth of the left-hand side compared to the linear growth dominated by the right-hand side leads to contradiction. In summary, when \( |p| \leq \frac{1}{s} \), behavior at infinity leads to contradiction.

Therefore, when \( p < 0 \), the problem (3.1) has no classical solution on \((x,t) \in \mathbb{R}^n_{1,+} \times \mathbb{R}\).

\section{Conclusion and Outlook}\label{sec4}

\quad \quad This paper has systematically investigated the existence, symmetry, monotonicity, and nonexistence of classical positive solutions to a class of elliptic and parabolic equations involving the fractional Laplacian with a nonlinearity of the form \(x_1 u^p + f(u)\). The analysis is carried out in both bounded domains (e.g., the unit ball \(B_1(0)\)) and the half-space \(\mathbb{R}^{n}_{1,+}\), covering three regimes: subcritical, supercritical, and negative exponent cases.

In the subcritical case (\(0 < p < \frac{n+2s}{n-2s}\)), we established the existence of bounded classical positive solutions via variational methods, specifically the mountain pass theorem and energy estimates. Using the method of moving planes, we further proved that these solutions exhibit symmetry with respect to the hyperplane \(\{x_1 = 0\}\) and monotonicity along the \(x_1\)-direction. These results hold under mild assumptions on the nonlinear term \(f(u)\), such as Lipschitz continuity and controlled growth.

In the supercritical case (\(p \geq \frac{n+2s}{n-2s}\)), we showed that no nontrivial classical positive solutions exist. This was achieved through a Pohozaev-type identity for the elliptic case and a contradiction argument near the boundary or at infinity for the parabolic case, leveraging the blow-up behavior of the nonlinear term.

In the negative exponent case (\(p < 0\)), we also proved nonexistence of classical solutions. For bounded domains, the singularity of \(u^p\) near the boundary leads to a contradiction with the regularity of the fractional Laplacian. For the half-space, both boundary behavior and growth at infinity were analyzed to rule out the existence of solutions.

The techniques developed here, including variational methods, the moving plane technique, energy estimates, blow-up analysis, and barrier arguments, are robust and can be extended to more general nonlocal operators and nonlinearities. The results not only deepen the theoretical understanding of fractional PDEs but also provide practical insights for modeling phenomena with long-range interactions and asymmetric nonlinearities.

Future work may include extending these results to more general domains, systems of equations, or other types of nonlocal operators, as well as exploring applications in physics, biology, and finance where such equations naturally arise.

\bibliographystyle{plain}
\bibliography{ref}

Lu Haipeng

School of Mathematics and Statistics

Northwestern Polytechnical University

Xi'an, Shaanxi, 710129, P. R. China

luhp2023@mail.nwpu.edu.cn

Yu Mei

School of Mathematics and Statistics

Northwestern Polytechnical University

Xi'an, Shaanxi, 710129, P. R. China

yumei@nwpu.edu.cn

\end{document}